\documentclass[oneside]{amsart}
\usepackage{epsfig,amsmath,latexsym,amssymb,amsfonts,amsthm,cite}
\usepackage{graphicx}
\usepackage[T1]{fontenc}
\usepackage{color,bm}
\usepackage{diagbox}
\graphicspath{{../img/}}
\newtheorem{theorem}{Theorem}[section]
\newtheorem{lemma}[theorem]{Lemma}
\newtheorem{proposition}[theorem]{Proposition}

\newtheorem{question}{Question}
\newtheorem{definition}{Definition}
\newtheorem{corollary}[theorem]{Corollary}
\newtheorem{main}{Theorem}

\def\F{\mathcal{F} }

\def\Q{\mathbb{Q} } 
\def\R{\mathbb{R} } 
\def\Z{\mathbb{Z} } 
\def\nbd{neighborhood } 
 
\def\R{\mathbb{R} } 
 
\def\Sv{\mathop{\mathrm{Sing}}(v)} 
\def\Pv{\mathop{\mathrm{Per}}(v)} 
\def\Cv{\mathop{\mathrm{Cl}}(v)}

\def\-{\ominus} 
\def\+{\oplus} 
\def\0{\circ}

\usepackage{color,bm}

\author{Tomoo Yokoyama}
\date{\today}
\address{Applied Mathematics and Physics Division, Gifu University, Yanagido 1-1, Gifu, 501-1193, Japan\\}
\email{tomoo@gifu-u.ac.jp}
\thanks{The author was partially supported by JSPS Grant Number 20K03583}
\makeatletter
\@namedef{subjclassname@2020}{%
  \textup{2020} Mathematics Subject Classification}
\makeatother

\title[Dependency of the positive and negative long-time behaviors]{Dependency of the positive and negative long-time behaviors of flows on surfaces} 



\begin{document}

\begin{abstract}
Long-time behavior is one of the most fundamental properties in dynamical systems. The limit behaviors of flows on surfaces are captured by the Poincar\'e-Bendixson theorem using the $\omega$-limit sets. This paper demonstrates that the positive and negative long-time behaviors are not independent. In fact, we show the dependence between the $\omega$-limit sets and the $\alpha$-limit sets of points of flows on surfaces, which partially generalizes the Poincar\'e-Bendixson theorem. Applying the dependency result to solve what kinds of the $\omega$-limit sets appear in the area-preserving (or, more generally, non-wandering) flows on compact surfaces, we show that the $\omega$-limit set of any non-closed orbit of such a flow on a compact surface is either a subset of singular points or a locally dense Q-set. Moreover, we show the wildness of surgeries to add totally disconnected singular points and the tameness of those to add finitely many singular points for flows on surfaces. 
\end{abstract}

\maketitle

\section{Introduction}

The Poincar\'e-Bendixson theorem is one of the most fundamental tools to capture the limit behaviors of orbits of flows and was applied to various phenomena (e.g. \cite{bhatia1966application,du2021traveling,koropecki2019poincare,hajek1968dynamical,pokrovskii2009corollary,roussarie2020topological,roussarie2021some,roussarie2021some02}). 
%
In \cite{birkhoff1927dynamical}, Birkhoff introduced the concepts of $\omega$-limit set and $\alpha$-limit set of a point. 
Using these concepts, one can describe the limit behaviors of orbits stated in the works of Poincar\'e and Bendixson in detail. 
Moreover, the Poincar\'e-Bendixson theorem was generalized for flows on surfaces in various ways  \cite{andronov1966qualitative,aranson1996introduction,buendia2017rem,buendia2019top,ciesielski1994poincare,demuner2009poincare,gardiner1985structure,gutierrez1986smoothing,Levitt1982foliation,lopez2004accumulation,lopez2007topological,markley1969poincare,marzougui1996,marzougui1998structure,nikolaev1999flows,schwartz1963generalization,vanderschoot2003local,yano1985asymptotic,yokoyama2021poincare}, and also for foliations \cite{levitt1987differentiability,plante1973generalization}, translation lines on the sphere \cite{koropecki2019poincare}, geodesics for a meromorphic connection on Riemann surfaces \cite{abate2016poincare,abate2011poincare}, group actions \cite{hounie1981minimal}, and semidynamical systems \cite{bonotto2008limit}. 
%
%
%
%
Applying a generalization of the Poincar\'e-Bendixson theorem, symmetric properties of long-time behaviors are studied \cite{yokoyama2023flows}.
 
Area-preserving flows on compact surfaces are one of the basic and classic examples of dynamical systems, also known as locally Hamiltonian flows or equivalently multi-valued Hamiltonian flows. 
The measurable properties of such flows are studied from various aspects  \cite{chaika2021singularity,conze2011cocycles,forni1997solutions,frkaczek2012ergodic,forni2002deviation,kanigowski2016ratner,kulaga2012self,ravotti2017quantitative,ulcigrai2011absence}. 
For instance, the study of area-preserving flows for their connection with solid-state physics and pseudo-periodic topology was initiated by Novikov \cite{novikov1982hamiltonian}. 
The orbits of such flows also arise in pseudo-periodic topology, as hyperplane sections of periodic manifolds (cf. \cite{arnol1991topological,zorich1999leaves}).

\subsection{Statements of main results}
In this paper, we discuss the dependency of the $\omega$-lilmit sets and the $\alpha$-limit sets of points and consider the following classification problem and the wildness and tameness of surgeries to add singular points of flows on surfaces.

\subsubsection{On the $\omega$-limit sets of points}

We show the following dependency of $\omega$-lilmit sets and the $\alpha$-limit sets, which supplements a generalization of the Poincar\'e-Bendixson theorem \cite[Theorem~A]{yokoyama2021poincare}. 

\begin{main}\label{lem:ld}
The $\omega$-limit set of any locally dense orbit for a flow with arbitrarily many singular points on a compact surface is either a nowhere dense subset of singular points or a locally dense Q-set.
\end{main}

The previous theorem says that the $\omega$-limit set and the $\alpha$-limit set of any point are not independent in general. For instance, the $\omega$-limit set of a point whose $\alpha$-limit set is a limit circuit is not a locally dense Q-set. 
The details of dependency are stated in \S 6 (see Tables~\ref{table:01} and \ref{table:02}). 
We apply the previous dependency result to the following question.

\begin{question}
What kinds of the $\omega$-limit sets of points do appear in the non-wandering, divergence-free, and Hamiltonian flows with arbitrarily many singular points on compact surfaces, respectively?
\end{question}

We answer that only nowhere dense subsets of singular points and locally dense Q-sets appear in such cases. 
Similar results hold for locally dense orbits of flows and non-closed orbits of gradient flows (see \S~\ref{sec:grad} for details). 

To describe precise statements of the results, we recall some concepts as follows. 
An orbit is {\bf closed} if it is singular or periodic, and it is {\bf locally dense} if its closure has a nonempty interior. 
An orbit $O$ is {\bf recurrent} if $O \subseteq \omega(O) \cup \alpha(O)$. 
A {\bf Q-set} is the closure of a non-closed recurrent orbit. 
Notice that a Q-set is also called a {\bf quasiminimal set}. 
We have the following dichotomy for any non-closed orbit of a non-wandering flow. 

\begin{main}\label{lem:nw}
The $\omega$-limit set of any non-closed orbit for a non-wandering flow with arbitrarily many singular points on a compact surface is either a nowhere dense subset of singular points or a locally dense Q-set.
\end{main}

We have the following observation for a Hamiltonian flow with arbitrarily many singular points on a (possibly non-compact) surface. 

\begin{main}\label{cor:ham}
The $\omega$-limit set of any non-closed orbit for a Hamiltonian flow with arbitrarily many singular points on a surface consists of singular points. 
In particular, there are no non-closed recurrent points of a Hamiltonian flow with arbitrarily many singular points on a surface. 
\end{main}


%
%

\subsection{Wildness and tameness of surgeries to add singular points}

We consider the following question for invariance of the Hamiltonian property by surgeries to add singular points.

\begin{question}\label{adding_sing}
Is the Hamiltonian property for flows on surfaces invariant under multiplying a bump function to a  Hamiltonian vector field?
\end{question}

We answer that the Hamiltonian property for flows on compact surfaces is invariant when only finitely many singular points are added, but that the property is not invariant even if totally disconnected singular points are added. 
%
More precisely, every non-zero vector field can be deformed into a vector field with wandering domains using a bump function with totally disconnected critical points as follows. 

\begin{main}\label{prop:45}
For any non-zero vector field $X$ on a compact surface $S$, there is a smooth function $f \colon S \to [0,1]$ such that 
\\
{\rm(1)} The vector field $fX$ is not non-wandering. 
\\
{\rm(2)} Every orbit of $X$ is a union of orbits of $fX$. 
\\
{\rm(3)} The critical point set $f^{-1}(0)$ is totally disconnected. 
\\
{\rm(4)} If the singular point set $\operatorname{Sing}(X)$ is totally disconnected, then so is $\operatorname{Sing}(fX)$. 
\end{main}

On the other hand, the Hamiltonian property of vector fields with finitely many singular points on compact surfaces is invariant under adding finitely many singular points. 

To state more precisely, we recall some concepts and a fact as follows. 
A nonempty subset $A$ of a topological space $X$ is a {\bf level set} of a function $f \colon X \to \R$ if there is a value $c \in \R$ such that $A = f^{-1}(c)$.  
Recall that, for any Hamiltonian vector field $X$ on a surface $S$ with the Hamiltonian $h \colon S \to \R$, the set of connected components of level sets of the restriction $h|_{S - \mathop{\mathrm{Sing}}(X)}$ is the set of orbits of the restriction $X|_{S - \mathop{\mathrm{Sing}}(X)}$ and is a codimension one foliation on the surface $S - \mathop{\mathrm{Sing}}(X)$. 
Therefore we introduce a pre-Hamiltonian flow as follows. 
\begin{definition}
A flow $v$ on an orientable surface is {\bf pre-Hamiltonian} if there is a continuous function $H \colon S \to \R$ such that the set of connected components of level sets of the restriction $H|_{S - \mathop{\mathrm{Sing}}(v)}$ is the set of orbits of the restriction $v|_{S - \Sv}$ and is a codimension one foliation on the surface $S - \mathop{\mathrm{Sing}}(v)$. 
\end{definition}

Then $H$ is called the {\bf pre-Hamiltonian} of $v$. 
We have the following equivalence under finiteness of singular points. 

\begin{main}\label{cor:characterization_ham_finite}
The following are equivalent for a flow $v$ with finitely many singular points on a compact surface $S$: 
\\
{\rm(1)} The flow $v$ is Hamiltonian. 
\\
{\rm(2)} The flow $v$ is pre-Hamiltonian. 
\end{main}

This implies that Hamiltonian property of flows with finitely many singular points on compact surfaces is invariant under adding finitely many singular points. 

The present paper consists of seven sections.
In the next section, as preliminaries, we introduce fundamental concepts.
In \S~3, classifications of the $\omega$-limit sets of points in the non-wandering and Hamiltonian flows on surfaces are demonstrated, and a remark is stated. 
In \S~4, the tameness of surgeries to add finitely many singular points is described. 
In \S~5, we demonstrate the wildness of surgeries to add totally disconnected singular points. 
In \S~6, we show the existence of non-forbidden pairs of $\alpha$-limit and $\omega$-limit sets of points. 
The final section remarks on Question~\ref{adding_sing}, time-reversal symmetric limit sets, and a construction like the Cherry flow box. 


\section{Preliminaries}

\subsection{Topological notion}
Denote by $\overline{A}$ the closure of a subset $A$ of a topological space and by $\partial A := \overline{A} - \mathrm{int}A$ the boundary of $A$, where $B - C$ is used instead of the set difference $B \setminus C$ when $B \subseteq C$.
%
A {\bf curve} is a continuous mapping $C: I \to X$ where $I$ is a non-degenerate connected subset of a circle $\mathbb{S}^1$.
A curve is {\bf simple} if it is injective.
We also denote by $C$ the image of a curve $C$.
A simple curve is a {\bf simple closed curve} if its domain is $\mathbb{S}^1$ (i.e. $I = \mathbb{S}^1$).
A simple closed curve is also called a {\bf loop}. 
An {\bf arc} is a simple curve whose domain is an interval. 
An {\bf orbit arc} is an arc contained in an orbit. 

By a {\bf surface}, we mean a paracompact two-dimensional manifold, that does not need to be orientable.
A subset of a compact surface $S$ is {\bf essential} if it is not null homotopic in $S^*$, where $S^*$ is the resulting closed surface from $S$ by collapsing all boundary components into singletons.

\subsection{Notion of dynamical systems}
A {\bf flow} is a continuous $\R$-action on a manifold.
From now on, we suppose that flows are on surfaces unless otherwise stated.
Let $v : \R \times S \to S$ be a flow on a surface $S$.
For $t \in \R$, define $v_t : S \to S$ by $v_t := v(t, \cdot )$.
For a point $x$ of $S$, we denote by $O(x)$ the orbit of $x$, $O^+(x)$ the positive orbit (i.e. $O^+(x) := \{ v_t(x) \mid t > 0 \}$),  $O^-(x)$ the negative orbit (i.e. $O^-(x) := \{ v_t(x) \mid t < 0 \}$).
A subset is {\bf invariant} (or saturated) if it is a union of orbits. 
A point $x$ of $S$ is {\bf singular} if $x = v_t(x)$ for any $t \in \R$ and is {\bf periodic} if there is a positive number $T > 0$ such that $x = v_T(x)$ and $x \neq v_t(x)$ for any $t \in (0, T)$.
A point is {\bf closed} if it is singular or periodic.
An orbit is singular (resp. periodic, closed) if it contains a singular (resp. periodic, closed) point. 
Denote by $\mathop{\mathrm{Sing}}(v)$ the set of singular points and by $\mathop{\mathrm{Per}}(v)$ (resp. $\mathop{\mathrm{Cl}}(v)$) the union of periodic (resp. closed) orbits. 

A point is {\bf wandering} if there are its neighborhood $U$
and a positive number $N$ such that $v_t(U) \cap U = \emptyset$ for any $t > N$.
A point is {\bf non-wandering} if it is not wandering (i.e. for any its neighborhood $U$ and for any positive number $N$, there is a number $t \in \mathbb{R}$ with $|t| > N$ such that $v_t(U) \cap U \neq \emptyset$).
Denote by $\Omega (v)$ the set of non-wandering points, called the non-wandering set. 
A flow is {\bf non-wandering} if any points are non-wandering.  

The $\omega$-limit (resp. $\alpha$-limit) set of a point $x$ is $\omega(x) := \bigcap_{n\in \mathbb{R}}\overline{\{v_t(x) \mid t > n\}}$ (resp.  $\alpha(x) := \bigcap_{n\in \mathbb{R}}\overline{\{v_t(x) \mid t < n\}}$). 
A {\bf separatrix} is a non-singular orbit whose $\alpha$-limit or $\omega$-limit set is a singular point.

A point $x$ is {\bf Poisson stable} (or strongly recurrent) if $x \in \omega(x) \cap \alpha(x)$. 
A point $x$ is {\bf positively} (resp. {\bf negatively}) {\bf recurrent} (or positively (resp. negatively) Poisson stable) if $x \in \omega(x)$ (resp. $x \in \alpha(x)$), and a point $x$ is {\bf recurrent} if $x \in \omega(x) \cup \alpha(x)$. 
Denote by $\mathrm{R}(v)$ the set of non-closed recurrent points.
An orbit is recurrent if it contains a recurrent point. 
Recall that a {\bf Q-set} is the closure of a non-closed recurrent orbit. 
Notice that a Q-set is also called a {\bf quasiminimal set}, and that a Q-set need not be orientable (see such examples \cite[Theorem~1]{gutierrez1978smooth}). 




An orbit is {\bf proper} if it is embedded, {\bf locally dense} if its closure has a nonempty interior, and {\bf exceptional} if it is neither proper nor locally dense. 
A point is proper (resp. locally dense) if its orbit is proper (resp. locally dense).
Denote by $\mathrm{LD}(v)$ (resp. $\mathrm{E}(v)$, $\mathrm{P}(v)$) the union of locally dense orbits (resp. exceptional orbits, non-closed proper orbits). 
By definitions, the union $\mathrm{P}(v)$ of non-closed proper orbits is the set of non-recurrent points, and that $\mathrm{R}(v) = \mathrm{LD}(v) \sqcup \mathrm{E}(v)$, where $\sqcup$ denotes a disjoint union (cf. \cite[\S 2.2.1]{yokoyama2021density}).
Moreover, we have a decomposition $S = \mathop{\mathrm{Cl}}(v) \sqcup \mathrm{P}(v) \sqcup \mathrm{R}(v) = \mathop{\mathrm{Sing}}(v) \sqcup \mathop{\mathrm{Per}}(v) \sqcup \mathrm{P}(v) \sqcup \mathrm{LD}(v) \sqcup \mathrm{E}(v)$ for a flow $v$ on a surface $S$. 
Every non-wandering flow on a compact surface has no exceptional orbits (i.e. $\mathrm{E}(v) = \emptyset$) because of \cite[Lemma~2.3]{yokoyama2016topological}.


\subsubsection{Quasi-circuits and quasi-Q-sets}
A closed connected invariant subset is a {\bf non-trivial quasi-circuit} if it is a boundary component of an open annulus, contains a non-recurrent orbit, and consists of non-recurrent orbits and singular points. 
A non-trivial quasi-circuit $\gamma$ is a {\bf quasi-semi-attracting quasi-circuit} if there is a point $x \in \mathbb{A}$ with $O^+(x) \subset \mathbb{A}$ such that $\omega(x)  = \gamma$. 
Moreover, a quasi-circuit is not a circuit in general. 


The transversality for a continuous flow can be defined using tangential spaces of surfaces, because each flow on a compact surface is topologically equivalent to a $C^1$-flow by Gutierrez's smoothing theorem~\cite{gutierrez1978structural}.
An $\omega$-limit (resp. $\alpha$-limit) set of a point is a {\bf quasi-Q-set} if it intersects an essential closed transversal infinitely many times. 
A quasi-Q-set is not Q-set in general, but a Q-set is a quasi-Q-set \cite[Lemma~3.8]{yokoyama2021poincare}.


\subsubsection{Types of singular points}
\begin{figure}
\begin{center}
\includegraphics[scale=0.325]{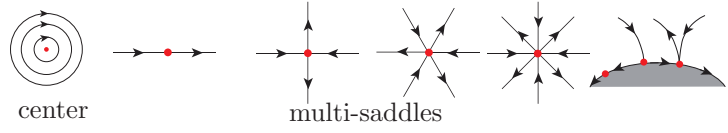}
\end{center}
\caption{A center and examples of multi-saddles}
\label{multi-saddles}
\end{figure}

A point $x$ is a {\bf center} if, for any its neighborhood $U$, there is an invariant open neighborhood $V \subset U$ of $x$ such that $U - \{ x \}$ is an open annulus that consists of periodic orbits, as in the left on Figure~\ref{multi-saddles}. 
A {\bf $\bm{\partial}$-$\bm{k}$-saddle} (resp. {\bf $\bm{k}$-saddle}) is an isolated singular point on (resp. outside of) $\partial S$ with exactly $(2k + 2)$-separatrices, counted with multiplicity as in Figure~\ref{multi-saddles}.
A {\bf multi-saddle} is a $k$-saddle or a $\partial$-$(k/2)$-saddle for some $k \in \mathbb{Z}_{\geq 0}$.
A $1$-saddle is topologically an ordinary saddle, and a $\partial$-$(1/2)$-saddle is topologically a $\partial$-saddle.

The union of multi-saddles and their separatrices is called the {\bf multi-saddle connection diagram}.
Any connected components of the multi-saddle connection diagram are called {\bf multi-saddle connections}.

\subsubsection{Hamiltonian flows}

A $C^r$ vector field $Y$ for any $r \in \Z_{\geq0}$ on an orientable surface $\Sigma$ is {\bf Hamiltonian} if there is a $C^{r+1}$ function $H \colon \Sigma \to \mathbb{R}$, called the {\bf Hamiltonian}, such that $dH= \omega(Y, \cdot )$ as a one-form, where $\omega$ is a volume form of $\Sigma$.
In other words, locally the Hamiltonian vector field $X$ is defined by $Y = (\partial H/ \partial x_2, - \partial H/ \partial x_1)$ for any local coordinate system $(x_1,x_2)$ of a point $p \in \Sigma$.
%
A flow is {\bf Hamiltonian} if it is topologically equivalent to a flow generated by a Hamiltonian vector field.

%
%

\subsubsection{Trivial flow boxes, flow boxes, periodic annuli, and transverse annuli}
\begin{figure}
\begin{center}
\includegraphics[scale=0.35]{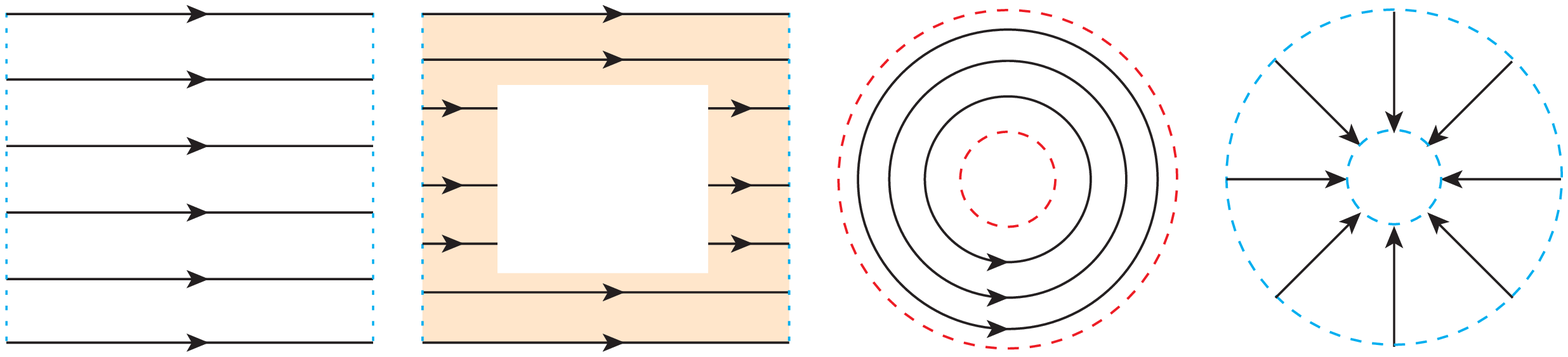}
\end{center}
\caption{Left, trivial flow box; second from left, flow box; second from right, periodic annulus; right, transverse annulus.}
\label{fig:local}
\end{figure}

A {\bf trivial flow box} is homeomorphic to $[0,1]^2$ each of whose orbit arcs correspond to $[0,1] \times \{t \}$ for some $t \in [0,1]$ as on the left of Figure~\ref{fig:local}. 
A closed subset $D$ is a {\bf flow box} with respect to a flow $v$ if there is a homeomorphism $h \colon D \to [0,1]^2$ 
 such that every connected component of the intersection $U \cap h^{-1}([0,1] \times \{t \})$ for any $t \in [0,1]$ is an orbit arc of $v$ as on second from left of Figure~\ref{fig:local}, where $U := h^{-1}([0,1]^2 - [1/4,3/4]^2)$ is a \nbd of the boundary $\partial D$ in $D$. 
A annulus $U$ is a {\bf periodic annulus} if there is a non-degenerate interval $I$ and there is a homeomorphism $h \colon U \to I \times \mathbb{S}^1$ such that the inverse image $h^{-1}(\{ t \} \times \mathbb{S}^1)$ for any $t \in I$ is a periodic orbit of $v$ 
as on second from right of Figure~\ref{fig:local}. 
A subset $U$ is an {\bf open transverse annulus} if there is a homeomorphism $h \colon U \to (0,1) \times \mathbb{S}^1$ such that the intersection $U \cap h^{-1}((0,1) \times \{ \theta \})$ for any $\theta \in \mathbb{S}^1$ is an orbit arc as on the right of Figure~\ref{fig:local}. 

\section{Proofs of main results and a remark on gradient cases}

We have the following observation. 

\begin{lemma}\label{lem:3.0}
Let $v$ be a flow with arbitrarily many singular points on a surface $S$ and $x$ a point in $S$. 
If $\omega(x) \subseteq \Sv$ {\rm(resp.} $\alpha(x) \subseteq \Sv${\rm)} then $\omega(x)$ {\rm(resp.} $\alpha(x)${\rm)} is a nowhere dense subset of singular points. 
\end{lemma}

\begin{proof}
If $x$ is singular, then $\omega(x) = \alpha (x) = \{ x \}$. 
Thus we may assume that $O(x)$ is not closed.  
Suppose that $\omega(x) \subseteq \Sv$. 
By $O(x) \cap \Sv = \emptyset$, we have $\omega(x) \subseteq \partial \Sv$. 
Since the boundary $\partial \Sv$ is closed and contains no interior, the $\omega$-limit set $\omega(x)$ is a nowhere dense subset of singular points. 
By the same argument, if $\alpha(x) \subseteq \Sv$, then $\alpha(x)$ is a nowhere dense subset of singular points. 
\end{proof}

We have the following statement, whose proof is analogous to the proof of \cite[Proposition~2.6]{yokoyama2016topological}. 

\begin{lemma}\label{lem:3.1}
Let $v$ be a flow with arbitrarily many singular points on a compact surface $S$. 
Then the following statements hold for a positively recurrent point $x$ in $\mathrm{LD}(v)$ and any point $y \in \overline{O(x)}$: 
\\
{\rm(1)} If $y$ is not positively recurrent, then $\omega(y) \subseteq \mathop{\mathrm{Sing}}(v)$. 
\\
{\rm(2)} If $y$ is not negatively recurrent, then $\alpha(y) \subseteq \mathop{\mathrm{Sing}}(v)$. 
\\
{\rm(3)} If $y \in \mathrm{P}(v)$, then $\alpha(y) \cup \omega(y) \subseteq \mathop{\mathrm{Sing}}(v)$. 
\end{lemma}

\begin{proof}
We may assume that the surface $S$ is connected. 
By taking a double covering of $S$ and the doubling of $S$ if necessary, we may assume that $S$ is closed and orientable. 
From \cite[Proposotion~2.3.11]{Hector1981foliation}, the restriction $v|_{S - \Sv}$ of the flow $v$ is transversally orientable as a foliation. 
%
Let $x$ be a positively recurrent point $x$ in $\mathrm{LD}(v)$ and a point $y \in \overline{O(x)}$. 
Then $O(y) \subset \omega(x) = \overline{O(x)}$. 
If $y \in \mathop{\mathrm{Sing}}(v)$, then $\{ y \} = \alpha(y) = \omega(y) \subseteq \mathop{\mathrm{Sing}}(v)$.
Thus we also may assume that $y \notin \mathop{\mathrm{Sing}}(v)$. 
\cite[Proposition~2.2 and Lemma~2.3]{yokoyama2016topological} imply that $\overline{O(x)} \cap (\mathop{\mathrm{Per}}(v) \sqcup \mathrm{E}(v)) = \emptyset$ and so that $y \in \overline{O(x)} \subseteq \mathop{\mathrm{Sing}}(v) \sqcup \mathrm{P}(v) \sqcup \mathrm{LD}(v)$. 
Therefore $y \in \mathrm{P}(v) \sqcup \mathrm{LD}(v)$. 
Since $\omega(y) = \omega(y_-)$ and $\alpha(y) = \alpha(y_-)$ for any point $y_- \in O^-(y)$, by replacing $y$ with a point in $O^-(y)$ if necessary, we may assume that $x \notin O^-(y) \sqcup \{ y \}$. 
Then $y \in (\mathrm{P}(v) \sqcup \mathrm{LD}(v)) - \{x\}$.


By the flow box theorem for a continuous flow on a compact surface (cf. \cite[Theorem~1.1, p.45]{aranson1996introduction}), there is a trivial flow box centered at the non-singular point $y$ as on the left in Figure~\ref{fig:flowbox_sing}.
\begin{figure}
\begin{center}
\includegraphics[scale=0.5]{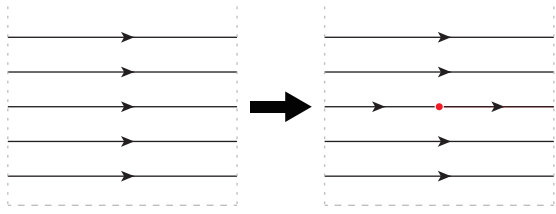}
\end{center}
\caption{Replacement of a trivial flow box on a disk into a flow fox with a singular point, which is a fake saddle.}
\label{fig:flowbox_sing}
\end{figure}
Replacing a trivial flow box on a disk into a flow fox with a singular point which is $y$ as on the right in Figure~\ref{fig:flowbox_sing}, denote by $v_1$ the resulting flow on the surface $S$. 
Notice that $O_v(z) = O_{v_1}(z)$ for any $z \notin O_v(y)$. 
Similarly, $O^+_v(z) = O^+_{v_1}(z)$ for any $z \notin O_v^-(y)$, and $O^-_v(z) = O^-_{v_1}(z)$ for any $z \notin O_v^+(y)$. 
For any $y_- \in O_v^-(y)$, we have $O_v^-(y) = O_{v_1}(y_-)$. 
Similarly, for any $y_+ \in O_v^+(y)$, we obtain $O_v^+(y) = O_{v_1}(y_+)$. 
By $x \notin O^-(y)$, we have that $O_v^+(x) = O_{v_1}^+(x)$ and so that $O_v^-(y) \subset \overline{O_v(x)} = \omega_v(x) = \omega_{v_1}(x) = \overline{O_{v_1}(x)}$. 
Then $x \in \mathrm{LD}(v_1)$. 
\cite[Proposition~2.2 and Lemma~2.3]{yokoyama2016topological} imply that $\overline{O_{v_1}(x)} \cap (\mathop{\mathrm{Per}}(v_1) \sqcup \mathrm{E}(v_1)) = \emptyset$ and so that $\overline{O_{v_1}(x)} \subseteq \mathop{\mathrm{Sing}}(v_1) \sqcup \mathrm{P}(v_1) \sqcup \mathrm{LD}(v_1)$. 
Moreover, $\overline{O_{v_1}(x)} \cap \mathrm{LD}(v_1) = \overline{O_{v_1}(x)} \setminus (\mathop{\mathrm{Sing}}(v_1) \sqcup \mathrm{P}(v_1)) = \{ z \in S \mid \overline{O_{v_1}(x)} = \overline{O_{v_1}(z)} \} = \{ z \in S \mid \overline{O_{v_1}(x)} = \overline{O_{v_1}^+(z)} \text{ or } \overline{O_{v_1}(x)} = \overline{O_{v_1}^-(z)} \}$ and $O_v^-(y) \subset \overline{O_{v_1}(x)} \setminus \mathop{\mathrm{Sing}}(v_1) \subseteq (\mathrm{P}(v_1) \sqcup \mathrm{LD}(v_1)) \cap \overline{\mathrm{LD}(v_1)}$. 

Suppose that $y$ is not positive recurrent (i.e. $y \notin \omega_v(y)$) with respect to $v$. 
If $x \in O_v^+(y) \sqcup \{ y \}$, then $y \in \overline{O_v(x)} =  \omega_v(x) = \omega_v(y)$, which contradicts that $y$ is not positively recurrent with respect to $v$. 
If $x \in \overline{O_v^+(y)} - (O_v^+(y) \sqcup \{ y \})$, then $x \in \omega(y)$ and so $y \in \overline{O_v(x)} =  \omega_v(x) \subseteq \omega_v(y)$, which contradicts that $y$ is not positively recurrent with respect to $v$. 
Thus $x \notin \overline{O_v^+(y)}$. 
Since $x \notin O_v^-(y) \sqcup \{ y \}$, we obtain $x \notin O_v(y)$.
Then $O_v(x) = O_{v_1}(x)$ and so $\overline{O_{v_1}(y_+)} = \overline{O_v^+(y)} \subsetneq \overline{O_{v}(x)} = \overline{O_{v_1}(x)} = \omega_{v_1}(x)$ for any $y_+ \in O_v^+(y)$.  
From $x \in \mathrm{LD}(v_1)$, applying \cite[Proposition~2.2]{yokoyama2016topological} to $v_1$, we have $O_v^+(y) = O_{v_1}(y_+) \subseteq \overline{O_{v_1}(x)} \setminus (\mathop{\mathrm{Sing}}(v_1) \sqcup \mathrm{LD}(v_1)) \subseteq \mathrm{P}(v_1)$ for any $y_+ \in O_v^+(y)$. 
Since $O_v^+(y) \subseteq \mathrm{P}(v_1)$, applying \cite[Theorem~3.1]{marzougui1996} to the restriction of $v_1$ to the complement $S - \mathop{\mathrm{Sing}}(v_1) = S - (\mathop{\mathrm{Sing}}(v) \sqcup \{y \})$ which is the orientable open surface equipped with the foliation $\{ O_{v_1}(z) \mid z \in S - \mathop{\mathrm{Sing}}(v_1) \}$, 
the orbit $O_{v_1}(y_+)$ of $v_1$ for any $y_+ \in O_v^+(y)$ is closed in the complement $S - \mathop{\mathrm{Sing}}(v_1)$ and so $\omega_v(y) = \omega_{v_1}(y_+) \subset \mathop{\mathrm{Sing}}(v_1) - \{ y \} = \mathop{\mathrm{Sing}}(v)$, because $y \notin \omega_v(y)$. 
This implies assertion {\rm(1)}. 

Suppose that $y$ is not negatively recurrent (i.e. $y \notin \alpha_v (y)$) with respect to $v$. 
If $x \in \overline{O_v^-(y)} - (O_v^-(y) \sqcup \{y \})$, then $x \in \alpha_v(y)$ and so $y \in \overline{O_v(x)} \subseteq \alpha_v(y)$, which contradicts that $y$ is not negatively recurrent with respect to $v$. 
Since $x \notin O_v^-(y) \sqcup \{ y \}$, we obtain $x \notin \overline{O_v^-(y)}$. 
%
Then $O^+_v(x) = O^+_{v_1}(x)$ and so that $O_{v_1}(y_-) = O_v^-(y) \subsetneq \overline{O_v(x)} = \omega_{v}(x) = \omega_{v_1}(x)$ for any $y_- \in O_v^-(y)$.  
From $x \in \mathrm{LD}(v_1)$, applying \cite[Proposition~2.2]{yokoyama2016topological} to $v_1$, we have $O_v^-(y) = O_{v_1}(y_-) \subseteq \overline{O_{v_1}(x)} \setminus (\mathop{\mathrm{Sing}}(v_1) \sqcup \mathrm{LD}(v_1)) \subseteq \mathrm{P}(v_1)$ for any $y_- \in O_v^-(y)$. 
Since $O_v^-(y) \subseteq \mathrm{P}(v_1)$, applying \cite[Theorem~3.1]{marzougui1996} to the restriction of $v_1$ to the complement $S - \mathop{\mathrm{Sing}}(v_1) = S - (\mathop{\mathrm{Sing}}(v) \sqcup \{y \})$ which is an orientable open surface equipped with an orientable and transversally orientable foliation $\{ O_{v_1}(z) \mid z \in S - \mathop{\mathrm{Sing}}(v_1) \}$, the orbit $O_{v_1}(y_-)$ of $v_1$ for any $y_- \in O_v^-(y)$ is closed in the complement $S - \mathop{\mathrm{Sing}}(v_1)$ and so $\alpha_v(y) = \alpha_{v_1}(y_-) \subset \mathop{\mathrm{Sing}}(v_1) - \{ y \} = \mathop{\mathrm{Sing}}(v)$, because $y \notin \alpha_v(y)$. 
This implies assertion {\rm(2)}. 

Suppose that $y \in \mathrm{P}(v)$. 
Then $y$ is not recurrent with respect $v$. 
Therefore assertions {\rm(1)} and {\rm(2)} imply $\alpha_v(y) \cup \omega_v(y) \subseteq \mathop{\mathrm{Sing}}(v)$.
\end{proof}

We demonstrate the main results for the $\omega$-limit sets as follows. 


%
%
%
\begin{proof}[Proof of Theorem~\ref{lem:ld}]
Let $v$ be a flow with arbitrarily many singular points on a compact surface $S$ and $y \in \mathrm{LD}(v)$. 
%
\cite[Theorem~VI]{cherry1937topological} implies that the Q-set $\overline{O(y)}$ contains infinitely many Poisson stable orbits $O$ with $\overline{O} = \overline{O(y)}$. 
Fix a Poisson stable point $x \in S$ with $\overline{O(x)} = \overline{O(y)}$. 
Then $x \in \mathrm{LD}(v)$ and $y \in \overline{O(x)} = \omega(x)$.  
If $y$ is not positively recurrent, then Lemma~\ref{lem:3.1} implies that $\omega(y) \subseteq \mathop{\mathrm{Sing}}(v)$ is a nowhere dense subset. 
Thus we may assume that $y$ is positively recurrent. 
Then $\omega(y) = \overline{O(y)}$ is a locally dense Q-set.

Therefore the $\omega$-limit set $\omega(x)$ is either a nowhere dense subset of singular points or a locally dense Q-set. 
\end{proof}



\begin{proof}[Proof of Theorem~\ref{lem:nw}]
Let $v$ be a non-wandering flow with arbitrarily many singular points on a compact surface $S$. 
By \cite[Lemma 2.4]{yokoyama2016topological}, we have that $S = \Cv \sqcup \mathrm{P}(v) \sqcup \mathrm{LD}(v) = \Sv \cup \overline{\Pv \sqcup \mathrm{LD}(v)}$. 
Fix a non-closed point $x \in S$. 

Suppose that $x \in \mathrm{LD}(v)$. 
If $x \in \omega(x)$, then $\omega(x) = \overline{O(x)}$ is a locally dense Q-set. 
Thus we may assume that $x \notin \omega(x)$. 
Then $x \in \alpha(x)$. 
By \cite[Theorem~VI]{cherry1937topological}, there is a Poisson stable point $x' \in \mathrm{LD}(v)$ with $\overline{O(x)} = \overline{O(x')}$. 
Since $x \in \overline{O(x')}$ is not positively recurrent, Lemma~\ref{lem:3.1} implies that the $\omega$-limit set $\omega(x) \subseteq \Sv$ is nowhere dense. 

Suppose that $x \notin \mathrm{LD}(v)$. 
Then $x \in S - (\Cv \sqcup \mathrm{LD}(v)) = \mathrm{P}(v)$. 
From \cite[Proposition~2.6]{yokoyama2016topological}, we have $\Sv \sqcup \mathrm{P}(v) = \{ y \in S \mid \omega(y) \cup \alpha(y) \subseteq \Sv \}$. 
Then the $\omega$-limit set $\omega(x)$ consists of singular points. 
By Lemma~\ref{lem:3.0}, the $\omega$-limit set $\omega(x) \subseteq \Sv$ is nowhere dense. 
\end{proof}

\subsection{Characterization of the $\omega$-limit set of any non-closed orbit for any non-wandering flow on a surface}
\label{sec:nw}

%
%
%

We have the following nonexistence of no non-closed recurrent orbits for Hamiltonian flows on (possibbly non-compact) surfaces. 


\begin{proof}[Proof of Theorem~\ref{cor:ham}]
Let $v$ be a Hamiltonian flow with arbitrarily many singular points on a surface $S$. 
By replacing $v$ with a flow which is topologically equivalent to $v$ and is generated by a Hamiltonian vector field, we may assume that there is a Hamiltonian $H \colon S \to \R$ generating $v$. 
Fix a non-closed point $x \in S$. 
Then the inverse image $H^{-1}(H(x))$ is a closed subset containing $\overline{O(x)}$. 
By the existence of the Hamiltonian, the point $x$ is proper and so $x \notin \omega(x) \cup \alpha(x)$. 
The invariance of $\omega(x)$ implies $O(x) \cap \omega(x) = \emptyset$. 

We claim that $\omega(x)$ consists of singular points. 
Indeed, assume that there is a non-singular point $y \in \omega(x)$. 
From $y \in \omega(x) \subseteq \overline{O(x)} \subseteq H^{-1}(H(x))$, we have $H(x) = H(y)$. 
%
By the flow box theorem for a continuous flow on a compact surface (cf. Theorem 1.1, p.45\cite{aranson1996introduction}), there is a closed disk $U$ which can be identified with $[-1,1] \times [-1,1]$ such that $y$ corresponds to the origin $(0,0)$ and the arc $C_t := [-1,1] \times \{ t \}$ for any $t \in [-1, 1]$ is contained in some orbit. 
By the definition of Hamiltonian vector field, for any $t \in [-1,1]$, there is a number $r_t \in \R$ with $H(C_t) = \{r_t\}$ and the function $r_{\cdot} \colon [-1,1] \to \R$ defined by $r_\cdot(t) := r_t$ is strictly increasing or decreasing. 
Moreover, we have $H(y) = r_0 \in \R$. 
On the other hand, since $O(x) \cap \omega(x) = \emptyset$, we obtain that $y \in \omega(x) \subseteq \overline{O(x)} - O(x)$ and so that the orbit $O^+(x)$ intersects $U - C_0 = U \setminus O(y)$.
Then there is a positive number $t^+ \in \R_{>0}$ such that $v_{t^+}(x) \in U - C_0$. 
This means that $H(v_{t^+}(x)) \neq r_0 = H(y) = H(x)$, which contradicts $H(x) = H(v_{t^+}(x))$. 
\end{proof}


\subsection{A remark on an analogous statement for gradient flows}\label{sec:grad}

Notice that a similar statement holds for a gradient flow with arbitrarily many singular points on a (possibly non-compact) surface. 
To state more precisely, recall the definition of gradient flows. 
A vector field $X$ on a Riemannian manifold $(M, g)$ is a {\bf smooth gradient vector} field if there is a $C^\infty$ function $f$ on $M$ with $g(X, \cdot) = d f$.
A flow is {\bf gradient} if it is topologically equivalent to a flow generated by a smooth gradient vector field.
Then we have the following statement. 

\begin{proposition}\label{lem:grad}
The $\omega$-limit set of any non-closed orbit for a gradient flow with arbitrarily many singular points on a surface consists of singular points.
\end{proposition}

\begin{proof}
Let $v$ be a gradient flow with arbitrarily many singular points on a surface $S$. 
By replacing $v$ with a flow which is topologically equivalent to $v$ and is generated by a gradient vector field, we may assume that there is a height function $h \colon S \to \R$ generating the gradient flow $v$. 
By the existence of the height function $h$, any orbits of $v$ are proper and non-periodic. 
%
%
Then $S = \Sv \sqcup \mathrm{P}(v)$. 
Fix a point $x \in  S - \Sv = \mathrm{P}(v)$. 

We claim that $\omega(x) \subseteq \Sv$.  
Indeed, 
assume that there is a point $y \in \omega(x) \cap \mathrm{P}(v)$. 
Then there are a closed transverse arc $I$ with $y \in I$ and an increasing sequence $(t_n)_{n \in \Z_{>0}}$ of $t_n \in \R_{>0}$ with $\lim_{n \to \infty} t_n = \infty$ such that $v_{t_n}(x) \in I$ and $y = \lim_{n \to \infty} v_{t_n}(x)$.   
By the existence of the height function $h \colon S \to \R$ of $v$, the value of $y$ is finite. 
Since $I$ is compact, the norms of the vector field $X_h := - \operatorname{grad} h$ generated by $h$ on some compact \nbd $U$ of $I$ are separated from zero (i.e. $\min_{z \in U}|X_h(z)| > 0$). 
This means that $h(y) = \lim_{n \to \infty} h(v_{t_n}(x)) = - \infty$, which contradicts $h(y) \in \R$.  
Thus $\omega(x) \subseteq S - \mathrm{P}(v) = \Sv$. 
%
\end{proof}

\section{On pre-Hamiltonian flows}

%

We have the following statements. 

\begin{lemma}\label{prop:characterization_ham_finite}
Every Hamiltonian flow on an orientable compact surface is pre-Hamiltonian. 
\end{lemma}

\begin{proof}
Let $v$ be a Hamiltonian flow on an orientable compact surface $S$. 
By definition of Hamiltonian flow, there is a Hamiltonian $h \colon S \to \R$ whose Hamiltonian vector field generates a flow $w$ which is topologically equivalent to $v$ via a topological conjugacy $k \colon S \to S$. 
By Theorem~\ref{cor:ham} and definition of Hamiltonian vector field, any orbits of the restriction $w|_{S - \mathop{\mathrm{Sing}}(w)}$ are closed in the surface $S -\mathop{\mathrm{Sing}}(w)$ and are connected components of level sets of the restriction $h|_{S - \mathop{\mathrm{Sing}}(w)}$. 
From the flow box theorem for a continuous flow on a compact surface, the set of orbits of $v|_{S - \mathop{\mathrm{Sing}}(w)}$ is a codimension one foliation on the surface $S - \mathop{\mathrm{Sing}}(w)$. 
Using the topological conjugacy $k \colon S \to S$ from $v$ to $w$, the composition $H := h \circ k$ is a desired continuous function $H \colon S \to \R$. 
\end{proof}

\begin{lemma}\label{lem:totally_ham_finite}
Let $v$ be a pre-Hamiltonian flow on an orientable compact surface $S$. 
If the image of $\mathop{\mathrm{Sing}}(v)$ by the pre-Hamiltonian is totally disconnected, then $v$ is non-wandering and $S = \Sv \cup \overline{\Pv}$. 
\end{lemma}

\begin{proof}
Let $v$ be a flow on an orientable compact surface $S$ and $H \colon S \to \R$ the pre-Hamiltonian of $v$. 
Since $\mathop{\mathrm{Sing}}(v)$ is closed and so compact, the image $H(\Sv)$ is compact and so closed. 

We claim that $S - H^{-1}(H(\Sv))$ consists of periodic orbits. 
Indeed, fix any point $x \in S$ with $H(x) \notin H(\Sv)$. 
Then the orbit $O(x)$ is a connected component of $H^{-1}(H(x))$ and so is a closed subset and so compact because $S$ is compact. 
Since $x$ is not a singular point, the orbit $O(x)$ is a periodic orbit. 

Fix any point $y \in H^{-1}(H(\Sv)) - \Sv$. 
Then the orbit $O(y)$ is a connected component of $H^{-1}(H(y)) \setminus \Sv$. 
The continuity of $H$ implies that every \nbd of $O(y)$ in $S - \Sv$ is not contained in $H^{-1}(H(y))$. 
By the totally disconnectivity of $H(\Sv)$, every \nbd of $O(y)$ in $S - \Sv$ intersects $H^{-1}(\R - H(\Sv)) = S - H^{-1}(H(\Sv)) \subseteq \Pv$. 
This means that $H^{-1}(H(\Sv)) -\Sv \subset \overline{\Pv}$ and so that $S = \Sv \cup \overline{\Pv}$. 
Therefore $v$ is non-wandering. 
\end{proof}

The total disconnectivity of the image $H(\Sv)$ of the singular point set by the pre-Hamiltonian is necessary and can not be replaced by one of the singular point set $\Sv$ in the previous lemma (see Corollary~\ref{cor:46}).  
To show Theorem~\ref{cor:characterization_ham_finite}, we recall the following concept.  
The {\bf extended orbit space $\bm{{S}/{v_{\mathrm{ex}}}}$} is a quotient space $S/\sim$ defined by $x \sim y$ if either $x$ and $y$ are contained in a multi-saddle connection or there is an orbit that contains $x$ and $y$ but is not contained in any multi-saddle connections. 
%
We prove Theorem~\ref{cor:characterization_ham_finite} as follows.


\begin{proof}[Proof of Theorem~\ref{cor:characterization_ham_finite}]
By Lemma~\ref{prop:characterization_ham_finite}, the assertion {\rm(1)} implies the assertion {\rm(2)}. 
Let $v$ be a pre-Hamiltonian flow on an orientable compact surface $S$, $H \colon S \to \R$ the pre-Hamiltonian, and $\F$ the foliation on $S - \Sv$ induced by $H|_{S - \Sv}$. 
By the Baire category theorem, any leaves of $\F$ have the empty interior. 
From Lemma~\ref{lem:totally_ham_finite}, the flow $v$ is non-wandering and $S = \Sv \cup \overline{\Pv}$. 
By \cite[Theorem 3]{cobo2010flows}, any singular points of $v$ are either centers or multi-saddles. 

We claim that there are no locally dense orbits. 
Indeed, assume that there is a locally dense orbit $O$. 
Then the closure $\overline{O}$ has a nonempty interior. 
Put $c := H(O) \in \R$. 
Since $\overline{O} \subseteq H^{-1}(H(O)) = H^{-1}(c)$, there is a connected component $L$ of the inverse image $H^{-1}(c)$ has a nonempty interior. 
By definition of pre-Hamiltonian, the connected component $L$ is a leaf of the codimension one foliation $\F$, which contradicts that every leaf has the empty interior. 

We claim that the extended orbit space ${S}/{v_{\mathrm{ex}}}$ is a directed graph without directed cycles.
Indeed, from \cite[Lemma 3.1]{yokoyama2021relations}, the multi-saddle connection diagram is the complement of the union of centers and $\Pv$, and the extended orbit space ${S}/{v_{\mathrm{ex}}}$ is a finite directed topological graph. 
Therefore any multi-saddle connections are connected components of $S - \Pv$. 
Since any level sets are closed, each multi-saddle connection is contained in some level set of $H$. 
This means that the quotient map $p_{\mathrm{ex}} \colon S \to {S}/{v_{\mathrm{ex}}}$ can be obtained by collapsing connected components of level sets of $H$ into singletons. 
Therefore there is an order-preserving continuous mapping $h \colon {S}/{v_{\mathrm{ex}}} \to \R$ with $H = h \circ p_{\mathrm{ex}}$. 
This implies that the directed graph ${S}/{v_{\mathrm{ex}}}$ has no directed cycles.

Since $v$ is a non-wandering flow without locally dense orbits whose extended orbit space ${S}/{v_{\mathrm{ex}}}$ is a directed graph without directed cycles, from \cite[Theorem~B]{yokoyama2021relations}, the flow $v$ is Hamiltonian.
\end{proof}

\section{Wildness of totally disconnected singular point sets}

We construct a non-trivial flow box with totally disconnected singular points that interrupts a flow. 

\subsection{Non-trivial flow box with totally disconnected singular points}\label{ex}


Recall that the {\bf Minkowski sum $\bm{A+B}$} is defined by $A+B := \{a+b \mid a \in A, b \in B\}$. 
Steinhaus shown that $\mathcal{C} + \mathcal{C} = [0,2]$, where $\mathcal{C}$ is the Cantor ternary set \cite{steinhaus1917new}. 
Using the Whitney theorem, we construct the following non-trivial flow box with totally disconnected singular points. 

\begin{lemma}\label{lem:tatally_disconnected}
Consider a vector field $X = \partial /\partial x^1 =  (1,0)$ on a closed square $D := [-1,2]^2$ and let $\mathbb{D}^2 := [0,1]^2 \subset D$. 
Then there is a smooth function $f \colon D \to [0,1]$ with $D - [-1/2, 3/2]^2 \subset f^{-1}(1)$ satisfying the following conditions: 
\\
{\rm(1)} The square $\mathbb{D}^2$ contains a wandering domain with respect to $fX$. 
\\
{\rm(2)} For any point $x \in \mathbb{D}^2$, we have either $O^+_{fX}(x) \subset \mathbb{D}^2$ or $O^-_{fX}(x) \subset \mathbb{D}^2$. 
\\
{\rm(3)} The intersection $\operatorname{Sing}(fX)$ is totally disconnected. 
\end{lemma}

\begin{proof}
Set $\mathcal{M} := \{ (x, x + y) \mid x,y \in \frac{1}{2} \mathcal{C} \} = \bigcup_{x \in \frac{1}{2} \mathcal{C}} \{ x\} \times \left( \{ x \} + \frac{1}{2}\mathcal{C} \right) \subset [0,1/2] \times [0,1] \subset \mathbb{D}^2$. 
By construction of $\mathcal{M}$, the subset $\mathcal{M}$ is a closed subset. 
By Whitney theorem \cite{Whitney1934analytic} (cf. \cite[Theorem~1.1.4]{Krantz1999geometry}), there is a $C^\infty$ function $f_0 \colon D \to [0,1]$ with $f_0^{-1}(0) = \mathcal{M}$. 
Using a bump function $\varphi \colon D \to [0,1]$ with $[0,1]^2 \subset \varphi^{-1}(0)$ and $D - [-1/2, 3/2]^2 \subset \varphi^{-1}(1)$, a $C^\infty$ function $f := f_0 (1-\varphi) + \varphi \colon D \to [0,1]$ satisfies $f^{-1}(0) = \mathcal{M}$ and $D - [-1/2, 3/2]^2 \subset f^{-1}(1)$. 
Then the vector field $fX$ is desired. 
\end{proof}

Using such a non-trivial flow box interrupting a flow, we show the deformation to create wandering domains.


\begin{proof}[Proof of Theorem~\ref{prop:45}]
By Lemma~\ref{lem:tatally_disconnected}, we can replace a trivial flow box $D$ with the non-trivial one as in  Lemma~\ref{lem:tatally_disconnected} by multiplying a smooth function $f \colon S \to [0,1]$ such that $f = 1$ outside of the flow box $D$, and that the restriction $f|_D$ is as in Lemma~\ref{lem:tatally_disconnected}. 
Thus the assertion holds. 
\end{proof}

We have the following statement, which shows the necessity of the total disconnectivity of the image of the singular point set by the pre-Hamiltonian in Lemma~\ref{lem:totally_ham_finite}. 

\begin{corollary}\label{cor:46}
There is a pre-Hamiltonian flow with wandering domains and totally disconnected singular points. 
\end{corollary}

\begin{proof}
Consider a Hamiltonian $h \colon \mathbb{S}^2 \to \R$ on the unit sphere $\mathbb{S}^2 = \{ (x,y,z) \in \R^3 \mid x^2+y^2+z^2 = 1 \}$ defined by $h(x,y,z) = z$. 
Replacing a trivial flow box for the Hamiltonian vector field $X$, by Theorem~\ref{prop:45}, there is a smooth function $f \colon \mathbb{S}^2 \to [0,1]$ such that the flow $v$ generated by the vector field $fX$ is not no-wandering and the singular point set $\Sv$ is totally disconnected. 
By construction, the function $h$ is the pre-Hamiltonian of $v$. 
\end{proof}

\section{Existence of pairs of the $\alpha$-limit and $\omega$-limit sets}

\begingroup
\renewcommand{\arraystretch}{1.4}
\begin{table}[htb]
  \begin{center}
\scalebox{0.9}{
    \begin{tabular}{|l|c|c|c|c|c|c|c|c|} 
    \hline
 \diagbox{$\alpha$-limit set}{$\omega$-limit set} & \begin{tabular}{c} Singular \\ point \end{tabular}  & \begin{tabular}{c} Limit\\ circuit \end{tabular}  & 
\begin{tabular}{c} Locally\\ dense \\Q-set \end{tabular} & \begin{tabular}{c} Transversely\\ Cantor \\ Q-set\end{tabular} \\  \hline  
Singular point & 1 & 2 &  4 & 5 \\  \hline  
Limit circuit & 2 & 3 & NO & 6\\  \hline  
Locally dense Q-set & 4 & NO & 7 &  NO \\  \hline  
\begin{tabular}{c} Transversely Cantor\\ Q-set\end{tabular} & 5 & 6 & NO  & 8 \\  \hline  
\end{tabular}
}
\end{center}
\caption{``NO'' represents the non-existence of pairs of the $\alpha$-limit and $\omega$-limit sets of any non-closed orbit for a flow with finitely many singular points on a compact surface.}\label{table:01}
\end{table}
\endgroup

\begingroup
\renewcommand{\arraystretch}{1.4}
\begin{table}[htb]
  \begin{center}
\scalebox{0.67}{
 \begin{tabular}{|l|c|c|c|c|c|c|c|c|} 
 \hline
 \diagbox{$\alpha$-limit set}{$\omega$-limit set } & \begin{tabular}{c} Nowhere \\dense subset\\ of \\ $\mathop{\mathrm{Sing}}(v)$  \end{tabular} & \begin{tabular}{c} Limit\\ cycle \end{tabular} & \begin{tabular}{c} Limit\\ quai-\\ circuit \end{tabular} & \begin{tabular}{c} Locally\\ dense \\Q-set \end{tabular} & \begin{tabular}{c} $\pitchfork$-ly\\ Cantor \\ Q-set\end{tabular} & \begin{tabular}{c} Quasi-Q-set in \\ $\Sv \sqcup \mathrm{P}(v)$\end{tabular}  \\  \hline  
\begin{tabular}{c} Nowhere dense subset\\ of $\mathop{\mathrm{Sing}}(v)$ \end{tabular} & $1'$ & 2 & 2 & 4 & 5 & $5'$ \\  \hline  
Limit cycle &2 & $3$ & $3'$ & NO & 6 & $6''$\\  \hline  
Limit quasi-circuit & 2& $3'$ & $3'$& NO & $6'$ & $6'''$ \\  \hline  
Locally dense Q-set & 4 & NO & NO & 7 &  NO  & NO  \\  \hline  
\begin{tabular}{c} $\pitchfork$-ly Cantor Q-set\end{tabular} & 5 & 6 & $6'$ & NO  & 8  & $8'$ \\  \hline  
\begin{tabular}{c} Quasi-Q-set in \\ $\Sv \sqcup \mathrm{P}(v)$ \end{tabular} &$5'$ &$6''$  &$6'''$ & NO  &  $8'$ & 9 \\  \hline  
\end{tabular}
}
\end{center}
\caption{``NO'' represents the non-existence of pairs of the $\alpha$-limit and $\omega$-limit sets of any non-closed orbit.}\label{table:02}
\end{table}
\endgroup

We show the existence of non-forbidden pairs of $\alpha$-limit and $\omega$-limit sets of points as follows. 

\begin{proposition}\label{prop:pair}
The pairs of the $\alpha$-limit and $\omega$-limit sets of points from $1$ to $9$ in Tables~\ref{table:01} and \ref{table:02} exist.
\end{proposition}

\begin{proof}
Cases $1$ and $1'$ occur in Morse flows on closed surfaces. 
Cases $2$ and $3$ occur in Morse-Smale flows on closed surfaces. 
Replacing a limit cycle with a limit quasi-circuit consisting of a singular point (or more generally, a simply connected domain of singular points) and one non-recurrent orbit, Case $3'$ occurs in the resulting flow. 
From an irrational rotation on a torus, replacing an orbit with the union of a singular point and two locally dense orbits, Case $4$ occurs in the resulting flow on the torus. 

From the suspension flow (i.e. Denjoy flow) of the Denjoy diffeomorphism, replacing a non-recurrent orbit with the union of a singular point and two non-recurrent orbits, Case $5$ occurs in the resulting flow on the torus. 

From a Denjoy flow and a Morse-Smale flows with limit cycles on a sphere, replacing an open flow box consisting contained in the union of non-recurrent orbits with a flow box as in Figure~\ref{fig:cherry01} for the Denjoy flow and the Morse-Smale flow respectively, by the time reversion for one of them, Case $6$ occurs in the resulting flow on the torus. 
\begin{figure}[t]
\begin{center}
\includegraphics[scale=0.25]{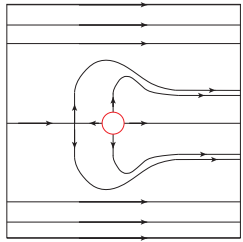}
\end{center}
\caption{A flow box on an annulus obtained from a Cherry flow box by replacing a sourc with a transverse boundary component.}
\label{fig:cherry01}
\end{figure}
Replacing a limit cycle with a limit quasi-circuit consisting of a singular point (or more generally, a simply connected domain of singular points) and one non-recurrent orbit, Case $6'$ occurs in the resulting flow. 

\cite[Theorem~VI]{cherry1937topological} implies that any Q-set contains infinitely many Poisson stable orbits and so Case $7$ occurs. 
In particular, Case $7$ occurs in any minimal flows on a torus. 

From a Denjoy flow, replacing an open flow box consisting contained in the union of non-recurrent orbits with a flow box as in Figure~\ref{fig:cherry01}, taking two copies of the resulting flow and by the time reversion for one of them, Case $8$ occurs in the resulting flow on the orientable closed surface $\Sigma_2$ whose genus is two. 

Let $\mathcal{C}$ be the Cantor minimal set of the Denjoy diffeomorphism and $\mathcal{M}$ the transversely Cantor minimal set of the suspension flow. 
As in \cite[Example~5.1]{yokoyama2021poincare}, replacing the Cantor set $\{0 \} \times \mathcal{C}$ with singular points by using a bump function, the subset $\mathcal{M}$ becomes a quasi-Q-set $\mathcal{M}'$consisting of singular points and separatrices of the resulting flow. 
%
%
Replacing a copy of $\mathcal{M}$ with $\mathcal{M}'$ on the torus in Case $5$ (resp. $6$, $6'$), Case $5'$ (resp. $6''$, $6'''$) occurs in the resulting flow. 
Replacing a copy of $\mathcal{M}$ with $\mathcal{M}'$ on $\Sigma_2$ in Case $8$, Case $8'$ occurs in the resulting flow. 
Similarly, replacing $\mathcal{M}$ with $\mathcal{M}'$ for the Denjoy flow on a torus, Case $9$ occurs in the resulting flow. 
\end{proof}

\section{Final remarks}

In this section, we state three remarks. 

\subsection{On Question~\ref{adding_sing}}

Note that the Hamiltonian property for flows on compact surfaces is invariant when only finitely many singular points are added. 
Indeed, consider a Hamiltonian flow $v$ on a compact surface. 
Adding a finitely many singular points to $v$, the resulting flow $v_1$ is pre-Hamiltonian. 
Theorem~\ref{cor:characterization_ham_finite} implies that the resulting flow $v_1$ is Hamiltonian. 
For instance, adding finitely many singular points to a closed periodic annulus, the resulting flow on the closed annulus is still Hamiltonian. 

\begin{figure}
\begin{center}
\includegraphics[scale=0.5]{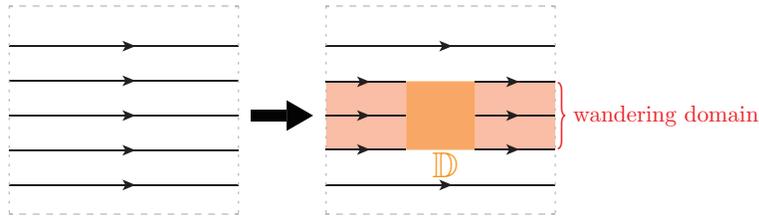}
\end{center}
\caption{Adding totally disconnected singular points to a trivial flow box.}
\label{fig:Cantor_sing}
\end{figure}

On the other hand, we can break the Hamiltonian property by adding totally disconnected singular points by Theorem~\ref{prop:45}. 
Indeed, replacing a trivial flow box to the flow box $D$ with the vector field $fX$ constructed in Lemma~\ref{lem:tatally_disconnected}, we obtain a wandering domain for the resulting flow. 
For instance, consider a closed annulus $A := \R/3\Z \times [-1,2]$ and immerse $D = [-1,2]^2$ into a subset of $A$ by the canonical projection $p \colon \R \times [-1,2] \to A$. 
Equip a vector field $X := (1,0)$ on $A$ and replace $X|_D$ with $fX|_D$, where $f \colon A \to [0,1]$ is the function constructed in Lemma~\ref{lem:tatally_disconnected}. 
The resulting flow $v$ satisfies $A = \Cv \sqcup \mathrm{P}(v)$, $\R/3\Z \times [0,1] = \Sv \sqcup \mathrm{P}(v)$, and $A - (\R/3\Z \times [0,1]) = \Pv$
  as in Figure~\ref{fig:Cantor_sing}. 

\subsection{Time-reversal symmetric limit sets}

In Table~\ref{table:01} and Table~\ref{table:02}, every case in the diagonal can be realized as time-reversal symmetric limit sets. 
Here an invariant subset $\mathcal{M}$ is a {\bf time-reversal symmetric limit set} \cite{yokoyama2023flows} if there is a point $x$ with $\mathcal{M} = \alpha(x) = \omega(x)$. 
More precisely, we have the following statement. 

\begin{proposition}
For any case in the diagonal in Table~\ref{table:01} and Table~\ref{table:02}, there is a flow on a compact surface with time-reversal symmetric limit sets in the case. 
\end{proposition}

\begin{proof}
Any flows with the singular points with homoclinic separatrices on compact surfaces are desired for Case~$1$ in Table~\ref{table:01}.

For Case $1'$ in Table~\ref{table:02}, we can construct flows with a time-reversal symmetric limit set which is a circle and consists of singular points. 
In fact, considering a Reeb component on a closed annulus as in Figure~\ref{fig:reeb_comp}, replacing the boundary into singular points, and gluing two boundary components, the resulting space is a torus and the resulting flow is continuous. 
By Gutierrez's smoothing theorem~\cite{gutierrez1978structural}, there is a smooth flow on a torus which is topologically equivalent to the resulting flow on the torus. 
\begin{figure}
\begin{center}
\includegraphics[scale=0.325]{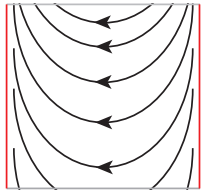}
\end{center}
\caption{Modification of a Reeb component.}
\label{fig:reeb_comp}
\end{figure} 
Notice that a circle consisting of singular points for a smooth flow on a closed disk is constructed in \cite[Example~3.2]{campos1997homeomorphisms}. 

A flow with a time-reversal symmetric limit cycle on a torus as in Figure~\ref{fig:reeb_nontwist} is desired for Case~$3$ in Table~\ref{table:01} and Table~\ref{table:02}. 
\begin{figure}
\begin{center}
\includegraphics[scale=0.3]{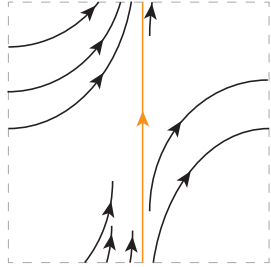}
\end{center}
\caption{A toral flow with one periodic orbit and non-recurrent orbits.}
\label{fig:reeb_nontwist}
\end{figure} 

For Case~$3'$ in Table~\ref{table:02}, there is a flow with a time-reversal symmetric limit set that has an invariant Wada-Lakes-like structure. 
In fact, there is a time-reversal symmetric limit set for a flow on a torus which is a limit quasi-circuit and is a complement of the union of an invariant open periodic center disk and an invariant open transverse annulus (see \cite[Example~4.3]{yokoyama2023flows}).

The minimal flows on a torus are desired for Case~$7$ in Table~\ref{table:01} and Table~\ref{table:02}. 
The Denjoy flow on a torus is desired for Case~$8$ in Table~\ref{table:01} and Table~\ref{table:02}. 
The flow in the proof of Proposition~\ref{prop:pair} in Case~9 is desired for Case~$9$ in Table~\ref{table:01}. 
\end{proof}

Notice that any locally dense $\omega$-limit sets in Case~7 need be time-reversal symmetric, and that there is a more complicated example for Case $1'$.
In fact, there is a non-wandering flow $v$ on an orientable closed surface consisting of singular points and locally dense orbits such that $\Sv$ is a ``double of lakes of Wada continuum'' and a time-reversal symmetric limit set \cite[Example~8.2]{yokoyama2023characterizations}. 
This construction implies that a nowhere dense subset of singular points in Theorem~\ref{lem:ld} and Theorem~\ref{lem:nw} can become a ``double of lakes of Wada continuum''.

\subsection{Construction like Cherry flow box}

\begin{figure}
\begin{center}
\includegraphics[scale=0.29]{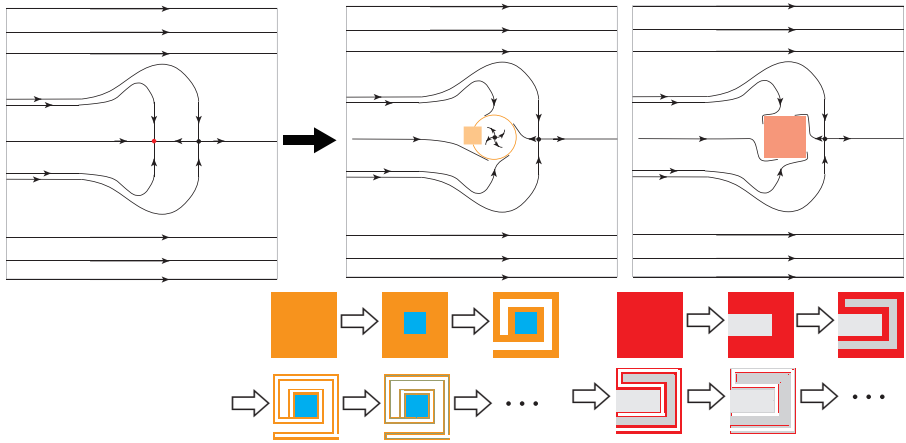}
\end{center}
\caption{Upper, Blowups of a Cherry flow box into flow boxes with $\omega$-limit sets each of which has an invariant Wada-Lakes-like structure; bottom middle and bottom right, the first few steps of a Wada-Lakes-like construction.}
\label{cherrybox_blowup}
\end{figure}

Using a limit quasi-circuit as above, we can construct a flow box with a $\omega$-limit set which is a limit quasi-circuit and has an invariant Wada-Lakes-like structure as on the middle in Figure~\ref{cherrybox_blowup} (see \cite[Example~4.3]{yokoyama2023flows} for details of the Wada-Lakes-like construction). 
Similarly, we can also construct flow box with a $\omega$-limit set which is a subset of the singular point set and is the Lakes of Wada as on the right in Figure~\ref{cherrybox_blowup}. 

\vspace{10pt}

\noindent
\textbf{Acknowledgements}

The author was partially supported by JSPS Grant Number 20K03583.

\bibliographystyle{plain}
\bibliography{../yt20220124}

\end{document}